\documentclass[letterpaper, 10 pt, conference]{ieeeconf}  

\IEEEoverridecommandlockouts                              

\usepackage{graphics} 
\usepackage{epsfig} 
\usepackage{mathptmx} 
\usepackage{times} 
\usepackage{amsmath} 
\usepackage{amssymb}  
\usepackage{amsfonts,xcolor}
\usepackage{cite}
\usepackage{array}

\usepackage{fourier-orns}

\title{\LARGE \bf
Transient Stability Analysis of Power Systems via Occupation Measures
}

\author{C\'edric Josz$^{1}$, Daniel K. Molzahn$^{2}$, Matteo Tacchi$^{3}$, and Somayeh Sojoudi$^{4}$
\thanks{*The research was funded by the DARPA grant D16AP00002, by the European Research Council (ERC) under the European Union's Horizon 2020 research and innovation program (grant agreement 666981 TAMING), and by the French company R\'{e}seaux de Transport d'\'{E}lectricit\'{e} (RTE).}
\thanks{$^{1}$C\'edric Josz is with the Department of Electrical Engineering and Computer Sciences, University of California, Berkeley, Etcheverry Hall, 2521 Hearst Ave, Berkeley, CA 94709.
        {\tt\small cedric.josz@gmail.com}}%
\thanks{$^{2}$Daniel K. Molzahn is with the Energy Systems Division, Argonne National Laboratory, Building 362, Room C381, 9700 South Cass Avenue, Lemont, IL 60439, USA.
        {\tt\small dmolzahn@anl.gov}}%
 \thanks{$^{3}$Matteo Tacchi is with the company R\'{e}seaux de Transport d'\'{E}lectricit\'{e}, 1 terrasse Bellini - Tour Initiale - TSA 41000 - 92919 Paris La D\'{e}fense Cedex, France ; as well as the LAAS-CNRS, 7 avenue du colonel Roche, 31031 Toulouse Cedex 4, France.
 {\tt\small tacchi@laas.fr}}%
\thanks{$^{4}$Somayeh Sojoudi is with Department of Electrical Engineering and Computer Sciences, University of California, Berkeley, 5114 Etcheverry Hall, 2521 Hearst Ave, Berkeley, CA 94709 and with the Tsinghua-Berkeley Shenzhen Institute.
 {\tt\small sojoudi@berkeley.edu}}%
}

\begin{document}

\maketitle
\thispagestyle{empty}
\pagestyle{empty}

\begin{abstract}

We propose the application of occupation measure theory to the classical problem of transient stability analysis for power systems. This enables the computation of certified inner and outer approximations for the region of attraction of a nominal operating point. In order to determine whether a post-disturbance point requires corrective actions to ensure stability, one would then simply need to check the sign of a polynomial evaluated at that point. Thus, computationally expensive dynamical simulations are only required for post-disturbance points in the region between the inner and outer approximations. We focus on the nonlinear swing equations but voltage dynamics could also be included. The proposed approach is formulated as a hierarchy of semidefinite programs stemming from an infinite-dimensional linear program in a measure space, with a natural dual sum-of-squares perspective. On the theoretical side, this paper lays the groundwork for exploiting the oscillatory structure of power systems by using Hermitian (instead of real) sums-of-squares and connects the proposed approach to recent results from algebraic geometry.

\end{abstract}


\section{Introduction}

The application of sum-of-squares (SOS) techniques to electric power systems dates back to 2000 in Parrilo's PhD thesis~\cite[Chapter 7.4]{parrilo-2000b}, where they are used for robust bifurcation analysis. More recently, there has been a growing interest in the power systems community regarding applications of SOS techniques and, in their dual form, moment relaxation hierarchies. In particular, these techniques are used to find global solutions to alternating current optimal power flow problems\cite{pscc2014,cedric_tps,ibm_opf,josz-phd}. 

The use of these techniques is justified when weaker relaxations~\cite{lavaei-low-2012} do not provide a global solution, but rather a strict lower bound~\cite{borden-demarco-lesieutre-molzahn-2011}. References~\cite{pscc2014,cedric_tps,ibm_opf,josz-phd} show that the Lasserre hierarchy of moment relaxations~\cite{lasserre-2001,parrilo-2003} can solve AC optimal power flow (ACOPF) problems for small power systems (with up to 10 buses) to global optimality using low orders of the hierarchy.
This is crucial since the Lasserre hierarchy becomes computationally expensive with increasing relaxation order.
By exploiting sparsity~\cite{mh_sparse_msdp,josz-phd,josz2018}, the Lasserre hierarchy can solve practical instances of ACOPF problems~\cite{josz-2016} with thousands of variables and constraints. This is achieved through a \emph{multi-ordered} Lasserre hierarchy~\cite{josz2018}. 

In this paper, we demonstrate that the problem of transient stability analysis (TSA) in power systems can be addressed using similar techniques. TSA considers the behavior of a power system following a major disturbance. The system must return to a stable condition and preserve synchronous operation after the switching of various devices and after faults. Electric power systems are growing in complexity due to increasing shares of renewable generation, increasing peak loads, and the expected wide-scale uses of demand response and energy storage. New tools are needed to benefit from high-performance computing and advances in sensing and communication equipment, such as phasor measurement units. Moreover, the control of power systems is complicated by phase-shifting transformers, HVDC lines, special protection schemes, etc. In this paper, we focus on uncontrolled dynamics as a first step towards certified estimations of the region of attraction (ROA) around a nominal operating point.

Similar to ACOPF problems, we find that TSA problems can be solved by convexifying the problem using measure theory, following the work of~\cite{korda2014} which admits a dual SOS perspective. To the best of our knowledge, SOS were first used to obtain estimates of the ROA of dynamical systems in~\cite[Chapter 7.3]{parrilo-2000b}. In the context of power systems, they were pioneered by the work~\cite{anghel2013}, which uses a Lyapunov approach (see~\cite{kundu2015acc,kundu2015cdc,kundu2015ecc,tacchi2018} for related works). The authors of~\cite{anghel2013} devise an expanding interior algorithm for estimating the ROA of the operating point. Their approach was recently improved in \cite{izumi2018} where an algorithm is devised that is simpler than the expanding interior, and includes convergence proofs, contrary to~\cite{anghel2013}. One could say that these previous approach are dual while the approach in this paper is primal. The key distinction is that the dual approach leads to sophisticated \textit{bilinear} matrix inequality conditions and relies on the choice of a shaping polynomial, while the primal approach results in a single semidefinite program with no additional data required besides the problem description and a hierarchy order. Moreover, the approach in \cite{izumi2018} only ensures convergence of the algorithm, but not necesssarily towards the the global optimum, while our primal approach based on \cite{korda2014} is endowed with a convergence in volume towards the actual ROA. We thus believe that it bears great potential for transmission systems operators, provided that sparsity may be exploited as in ACOPF problems.

We next summarize some recent work on power systems TSA. Wang~\text{et al.}~\cite{wang2017} propose TSAs using a hybrid direct-time-domain method and a partial energy function. The analysis of the power system is reduced to several pairs of ``coupled'' machines with large rotor speed differences. Owusu-Mireku and Chiang~\cite{mireku2017} propose an energy-based method for the TSA after a power system transmission switching event. Their method determines a relevant controlling unstable equilibrium point for a switching event and then uses an energy margin to assess stability. Dasgupta and Vaidya~\cite{dasgupta2018} develop a methodology for finite-time rotor TSA. The authors draw on the theory of normal hyperbolic surfaces in order to bring new insights to existing techniques for finite-time stability.
All these contributions are confirmed numerically on relevant test cases, such as those in~\cite{josz-2016}.


This paper is organized as follows. Section~\ref{sec:PROBLEM FORMULATION} formulates the TSA problem. 
Section~\ref{sec:PROPOSED APPROACH: OCCUPATION MEASURE} presents the proposed occupation-measure-based method as well as some foundational theoretical results. Section~\ref{sec:NUMERICAL EXPERIMENT} describes numerical experiments conducted to show the practical relevance of the proposed method and gives future research directions regarding computational tractability.

\section{Problem Formulation}

\label{sec:PROBLEM FORMULATION}

\subsection{Transient stability of power systems}
\label{subsec:overview}

Consider a power system composed of $n$ synchronous generators with respective complex voltages $v_1,\hdots,v_n$. We assume, as it is common in the literature, that the voltage magnitudes $|v_1|,\hdots,|v_n|$ are fixed during the transient period, while the phase angles $\theta_1,\hdots,\theta_n$ are variable (compared to the rotating frame) with respective angular speeds $\omega_1,\hdots,\omega_n$. In addition, the loads in the network are considered to be constant and passive impedances. After a fault occurs, the phases will satisfy the following set of differential equations:
\begin{equation} 
\label{eq:dyn_syst}
\left\{
\begin{array}{rcl} 
\dot{\theta}_k  & = & \omega_k, \\[0.25em]
\dot{\omega}_k & = & - \lambda_k \omega_k + \frac{1}{M_k}\left(P_k^\text{mec} - P_k^\text{elec}(\theta_1,\hdots,\theta_n) \right),
\end{array}
\right.
\end{equation}
where $P_k^\text{mec}$ is the (fixed) mechanical power input at bus~$k$ and $P_k^\text{elec}(\theta_1,\hdots,\theta_n)$ is the electrical power output of each generator~$k$ with value given by
\begin{equation} G_{kk} |V_k|^2 + \sum\limits_{l \neq k} |v_k| |v_l| \left\{ B_{kl} \sin(\theta_k - \theta_l) + G_{kl} \cos ( \theta_k - \theta_l ) \right\}.
\end{equation}
The quantities $B_{kl}$ and $G_{kl}$ denote the line susceptances and conductances, and $M_k$ denotes the generator inertia constant. The constant $\lambda_k$ is related to the damping coefficient of each generator.

We assume that there exists an equilibrium to these equations, i.e., values of $\theta^{\text{eq}}$ that satisfy

\begin{equation}
\label{eq:equilibrium}
P_k^\text{mec} = P_k^\text{elec}(\theta_1^\text{eq},\hdots,\theta_n^\text{eq}) , ~~~ k = 1,\hdots,n.
\end{equation}
In other words, $\theta^{\text{eq}}$ corresponds to a steady-state operating point of an AC transmission system. As usual, we choose one bus, denoted by subscript ``ref'', to serve as the reference bus, with $\theta^\text{eq}_{\text{ref}} = 0$ (often referred to as slack bus). Indeed, the equations are invariant up to a phase shift.  Although the focus of the paper is on frequency analysis, the results apply to a more comprehensive model coupled with voltage dynamics. The details are omitted for brevity.


The TSAs described in this paper rely on polynomial reformulations of the dynamical system model~\eqref{eq:dyn_syst}--\eqref{eq:equilibrium}. To illustrate these reformulations, we use the three-bus example from Chiang~\text{et al.}~\cite{chiang2011}, which is composed of three synchronous machines connected in a cycle. Since the third bus sets the reference angle (i.e., $\theta_3 = 0$), we only need two phase angle variables, $\theta_1,\,\theta_2$, and two rotor speed variables, $\omega_1,\,\omega_2$, to describe the dynamics:
$$
\left\{
\begin{array}{rcl}
\dot{\theta}_k &=& \omega_k, \qquad k=1,2, \\[0.25em]
\dot{\omega}_1 &=& - \sin(\theta_1) - 0.5 \sin (\theta_1-\theta_2) - 0.4 \omega_2,\\[0.25em]
\dot{\omega}_2 &=& - 0.5 \sin(\theta_2) - 0.5 \sin (\theta_2-\theta_1) - 0.5 \omega_2 + 0.05.
\end{array}
\right.
$$
A stable equilibrium is given by $(\theta_1^\text{eq},\,\theta_2^\text{eq}) = (0.02,\, 0.06)$. Following \cite{anghel2013}, the coordinates can be shifted so that $(\theta_1^\text{eq},\theta_2^\text{eq}) = (0.00,\, 0.00)$ is a stable equilibrium. 
This dynamical system can in turn be formulated as a polynomial differential algebraic system, as suggested by Anghel~\text{et al.}~\cite{anghel2013}. To that end, we introduce auxiliary variables 
\begin{equation}
\label{eq:varchange}
s_k := \sin(\theta_k) ~~~ \text{and} ~~~ c_k := 1-\cos(\theta_k), ~~~ k=1,2
\end{equation}
The reformulated dynamical system is
$$
\left\{
\begin{array}{cclr}
\dot{\omega}_1 \!\!\!&=&\!\!\! 0.4996s_2 - 0.4\omega_1 - 1.4994s_1 - 0.02 c_2 & \hspace{-0.9em}+ 0.02s_1s_2 \\[0.25em]
\!\!\!& &\!\!\! + 0.4996s_1c_2 - 0.4996 c_1 s_2 + 0.02c_1c_2, \\[0.25em]
\dot{\omega}_2 \!\!\! &=&\!\!\! 0.4996 s_1 + 0.02 c_1 - 09986 s_2 + 0.05 c_2 & \hspace{-1em} - 0.5\omega_2 \quad\hphantom{.} \\[0.25em]
\!\!\!& &\!\!\! -0.02 s_1s_2 - 0.4996 s_1c_2 + 0.4996 c_1s_2 & \hspace{-0.9em} - 0.02 c_1c_2, \\[0.25em]
\dot{s}_k \!\!\!&=&\!\!\! (1 - c_k)\omega_k & k=1,2,\\[0.25em]
\dot{c}_k \!\!\!&=&\!\!\! s_k \omega_k & k=1,2,\\[0.25em]
0 \!\!\!&=&\!\!\! s_k^2 + c_k^2 - 2.0 c_k & k=1,2,\\[0.25em]
\end{array}
\right.
$$
Section~\ref{sec:NUMERICAL EXPERIMENT} will show that one can actually avoid increasing the number of variables and immediately obtain an algebraic differential system of equations in \textit{complex-valued quantities}. 

\subsection{Region of attraction}

Consider the basic semi-algebraic set 
\begin{equation}
X := \{ ~ x \in \mathbb{R}^n ~|~ g_i(x) \geqslant 0, ~ i = 1,\hdots,n_X \}
\end{equation}
where $g_1,\hdots, g_{n_X}$ are polynomials such that $X$ is compact, as well as the differential algebraic system
\begin{equation}
\label{eq:dynamics}
\begin{cases}\dot x(t) = f(x(t)),\\ g_0(x(t))=0\end{cases}~~~x(t) \in X,~~~ \forall t \in [ 0 , T ],
\end{equation}
where $x(\cdot): [ 0 , T ] \longrightarrow \mathbb{R}^n$, $f \in \mathbb{R}[x]^n$, $T> 0$ and $g_0\in\mathbb{R}[x]$.

In addition, we ask that the final state $x(T)$ belongs to another semi-algebraic set $X_T \subset X$, for example, a Euclidian ball with a small radius $\epsilon>0$ centered at the equilibrium.


The region of attraction (ROA) $X_0$ is the set of initial conditions for which there exists an admissible trajectory:
\begin{align}
    \nonumber X_0 := & \left\{ ~ x_0 \in X ~\vert~ \exists ~ x(\cdot\vert x_0) \text{ solution to \eqref{eq:dynamics} on $[0,T]$ s.t. } \right.\\
    \nonumber & \left. \qquad x(0\vert x_0) = x_0, \text{ and } x(T\vert x_0) \in X_T \right\}.
\end{align}
The remainder of this paper describes approaches for computing inner and outer approximations to the ROA $X_0$.

\section{Approximation of the Region of Attraction via Occupation Measures}

\label{sec:PROPOSED APPROACH: OCCUPATION MEASURE}

In this section, we explain the general approach proposed by Henrion and Korda~\cite{korda2014,korda2016}. Their idea is to provide a convex formulation of polynomial ODEs using the notion of \emph{occupation measures} (OM)~\cite{vinter1993}, which quantify the time spent by the trajectory of the state in a set $B \subset X$:
\begin{equation}
\mu(A \times B | x_0 ) := \int_0^T I_{A \times B}(t,x(t|x_0))\, dt
\end{equation}
where $A \subset [0,T]$ and $I$ is the indicator function. Importantly, such a $\mu$ satisfies, for any measurable function \mbox{$\varphi:X\rightarrow \mathbb{R}$,}
\begin{equation}
\label{eq:occupation}
\int_0^T \varphi(t,x(t|x_0))\, dt = \int_{[0,T] \times X} \varphi(t,x)\,d \mu (t,x|x_0).
\end{equation}
Next, define the operator $\mathcal{L}: C^1( [ 0,T ] \times X) \rightarrow C([0,T] \times X)$
\begin{equation}
\label{def:ruelle}
v \longmapsto \mathcal{L}v := \frac{\partial v}{\partial t} + \sum_{i=1}^n \frac{\partial v}{\partial x_i} f_i(t,x) = \frac{\partial v}{\partial t} + \mathrm{grad}~ v \cdot f.
\end{equation}
Then, for any $v \in C^1([0,T]\times X,\mathbb{R})$, \eqref{eq:occupation} and \eqref{def:ruelle} yield
\begin{equation}
\label{eq:int}
v(T,x(T|x_0)) = v(0,x_0) + \int_{[0,T] \times X} \mathcal{L} v(t,x)\, d\mu(t,x|x_0).
\end{equation}
If instead of an initial point $x_0$, we consider a probability distribution $\mu_0$ supported on the feasible set $X$, one may define the \textit{average occupation measure}
\begin{align}
\mu (A \times B) & := \int_X \mu (A \times B | x_0)\, d\mu_0(x_0), \\
\mu_T(B) & := \int_X I_B (x(T|x_0))\, d\mu_0(x_0).
\end{align}
Integrating \eqref{eq:int} with respect to $\mu_0$, we obtain that 
\begin{align}\nonumber
\int_{X} v(T,x)\, d\mu_T(x) = & \int_X v(0,x)\, d\mu_0(x)  \\[-0.4em]
\label{eq:Liouville} & \quad + \int_{[0,T] \times X} \mathcal{L} v(t,x)\, d\mu(t,x) .
\end{align}
Using distributional derivatives, one can interpret the above equation as Liouville's PDE. Finding the ROA is then formulated as the following optimization problem:\vspace{-0.7em}
\begin{eqnarray}
\label{eq:optliouville}
p^\ast = & \sup & \mu_0 ( X ) \\
\nonumber & \text{s.t.} & \text{Liouville equation~\eqref{eq:Liouville}}, \\
&& \mu_0 + \hat{\mu}_0 = \lambda, \label{eq:slack}\\
\nonumber && \mu \geqslant 0 , ~ \mu_0 \geqslant 0 , ~ \mu_T \geqslant 0, ~ \hat{\mu}_0 \geqslant 0, \\
\nonumber && \mathrm{spt} (\mu) \subset [0,T] \times X, ~ \mathrm{spt} (\hat{\mu}_0) \subset X \\
\nonumber && \mathrm{spt}(\mu_0) \subset X , ~ \mathrm{spt}(\mu_T) \subset X_T.
\end{eqnarray}
where $\lambda$ denotes the Lebesgue measure on $X$ and $\mathrm{spt}$ denotes the support of a  measure. Equations \eqref{eq:Liouville} and \eqref{eq:slack} induce a linear relationship between the four measures. The optimal value of this infinite dimension linear program is equal to the volume of the ROA \cite[Theorem~1]{korda2014}. Importantly, the supremum is attained and the optimal solution is such that $\mu_0^\ast$ is the restriction of the Lebesgue measure to the ROA.


In his seminal article \cite{lasserre-2001}, Lasserre showed that such infinite-dimensional linear program on measures $\mu$ can be approximated by a hierarchy of finite-dimensional semidefinite programs on vectors of moments $y_\alpha = \int x^\alpha d\mu(x)$, \mbox{$\vert \alpha \vert \leqslant 2k$} \cite{lasserre-2010}. These hierarchies have the remarkable property of yielding upper bounds $p^\ast_k$ of the infinite-dimensional optimal value $p^\ast$ such that $p^\ast_k \underset{k \rightarrow \infty}{\searrow} p^\ast$.

There exists a dual perspective to the approach:
\begin{equation}
\begin{array}{rll}
d^\ast = \inf & \int_X w(x)\, d\lambda(x) \\[0.3em]
\text{s.t.} & \mathcal{L} v(t,x) \leqslant 0 , & \forall (t,x) \in [0,T] \times X, \\[0.3em]
& w(x) \geqslant v(0,x) + 1, & \forall x \in X, \\[0.3em]
& v(T,x) \geqslant 0 , & \forall x \in X_T, \\[0.3em]
& w(x) \geqslant 0 , &  \forall x \in X.
\end{array}
\end{equation}
The constraint $\mathcal{L} v(t,x) \leqslant 0$ implies that $v$ is non-increasing along the trajectories, and thus $v(0,x) \geqslant 0$ on $X_0$ due to the constraint $v(T,x) \geqslant 0$ on $X_T$. As a byproduct, we also have that $w(x) \geqslant 1$ on $X_0$. A nice property about the previous optimization problems is that there is no duality gap \cite[Theorem~2]{korda2014}.

This dual perspective naturally admits a SOS reformulation:
\begin{equation}
\begin{array}{rl}
\inf & \textbf{w}^\intercal\, h \\[0.3em]
\text{s.t.} & - \mathcal{L}v_k(t,x) = p(t,x) +q_0(t,x) t(T-t) \\[0.3em]
& \qquad\qquad\qquad + \sum_{i=1}^{n_X} q_i(t,x) g_i^X(x), \\[0.3em]
& w_k(x) - v_k(0,x) - 1 = p_0(x) + \sum_{i=1}^{n_X} q_{0_i}(x) g_i^X(x), \\[0.3em]
& v_k(T,x) = p_T(x) + \sum_{i=1}^{n_T} q_{T_i}(x) g_i^{X_T}(x), \\[0.3em]
& w_k(x) = s_0(x) + \sum_{i=1}^{n_X} s_{0_i}(x)g_i^X(x).
\end{array}
\end{equation}
where $h$ is the vector of $\lambda$'s moments, and $\textbf{w}$ is the vector of coefficients of $w_k(x)$ in the moments basis. The optimization variables include polynomials $v_k(t,x)$ and $w_k(x)$ of degree at most $2k$ as well as the SOS polynomials $p(x)$, $q_i(x)$, $p_0(x)$, $p_T(x)$, $q_{0_i}(x)$, $q_{T_i}(x)$, $s_0(x)$, and $s_{0_i}(x)$ with appropriate degrees that can be deduced from the constraints in the optimization problem. Again, there is no duality gap between the truncated problems at every order of the hierarchy~\cite[Theorem~4]{korda2014}.

An outer approximation to the ROA is then given by 
\begin{equation}
\boxed{\tilde{X}_0 := \{ ~ x \in \mathbb{R}^n ~|~ v_k(0,x) \geqslant 0 ~ \}}
\end{equation}
which converges in volume towards the ROA as the order~$k$ increases to infinity \cite[Theorem~6]{korda2014}. As with the Lasserre hierarchy or the Lyapunov approach via SOS, the computational burden increases sharply as the order $k$ increases.

A particularity of the OM approach is that the state set $X$ should have an interior point such that the computed volumes are non-zero. Hence, constraints $g_0(x(t)) = 0$ in \eqref{eq:dynamics} derived from our change of variable may be troublesome, since the manifold $M := \{ x \in X ~\vert ~ g_0(x) = 0 \}$ has no interior point. A simple method to address this issue consists in ignoring the equality constraints when computing the ROA approximation $\tilde{X}_0$, and then consider $\tilde{X}_0 \cap M$ as the desired ROA estimation. Such a method does not work with any arbitrary  equality constraints. However, in the case of constraints derived from a change of variable, this approach is valid due to the fact that the vector field $f$ then satisfies $(\mathrm{grad}~g_0) \cdot f \equiv 0.$ Thus, the dynamics are tangent to $M$, which means that any trajectory starting in $M$ will remain in $M$, which is exactly the constraint $g_0(x(t))=0, ~ \forall ~ t \in [0,T]$.

To the best of our knowledge, this is the first time that algebraic equality constraints derived from a change of variable are addressed within the OM approach. This facilitates the novel application of OM theory to non-polynomial systems.

We conclude this section by briefly discussing the approach for computing inner approximations. The machinery for inner approximations is very similar to the outer approximation approach discussed above. The key distinction is that the inner approximations consider an outer approximation to the complement of the ROA, $X_0^c := X \setminus X_0$. See~\cite{korda2013} for further details.

   \begin{figure}[t]
      \centering
      \includegraphics[width=9cm]{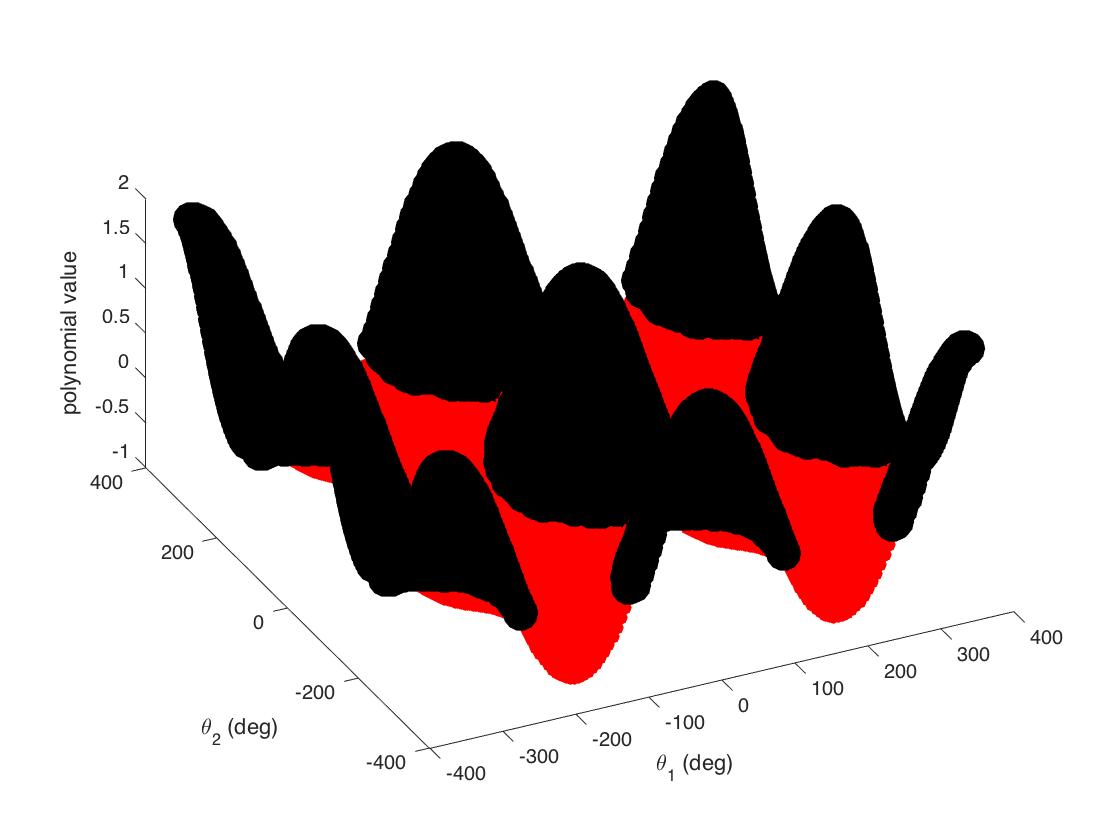}
      \caption{The polynomial for the three-bus system whose zero level set, which is indicated by the back region, provides an outer approximation to the ROA. The projection shown is for $(\omega_1,\,\omega_2) = (0,\,0)$.}
      \label{fig:level}
      \vspace*{-1em}
   \end{figure}
   
\vspace{-0.7em}
\section{Case study}

\label{sec:NUMERICAL EXPERIMENT}

For our numerical experiments, we use MATLAB R2015b, YALMIP~\cite{yalmip}, SeDuMi~1.3~\cite{sturm-1999}, and the ``ROA'' code of Henrion and Korda~\cite{korda2014} to apply OM theory to the three-bus example from~\cite{chiang2011} that is described in Section~\ref{subsec:overview}.

We note that practical power system analyses require the ability to address significantly larger problems than the test case considered in this paper. However, constructing certified approximations for the ROA leads to difficult computational challenges. Similar to the demonstrations of previous algorithms~\cite{anghel2013},\cite{izumi2018}, this paper focuses on a small system as an initial step towards practical applications. Future work that exploits network sparsity and other problem structures will be crucial for scalability. Decomposition approaches may also prove valuable~\cite{kundu2015ecc,kundu2015cdc,kundu2015acc}.

   \begin{figure}[t]
      \centering
      \includegraphics[width=9.25cm]{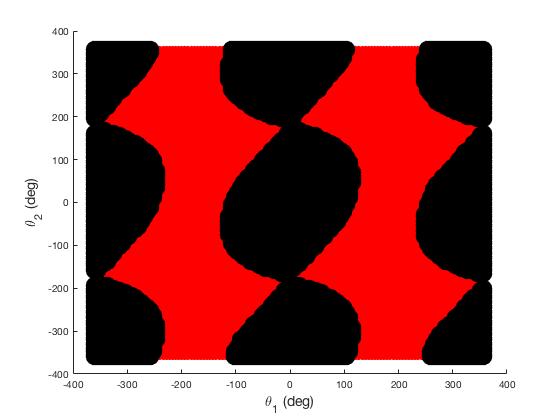}
      \caption{An outer approximation of the ROA is indicated by the back region. The projection shown is for $(\omega_1,\,\omega_2) = (0,\,0)$.}
      \label{fig:outer}
   \end{figure}

   \begin{figure}[t]
      \centering
      \includegraphics[width=9.25cm]{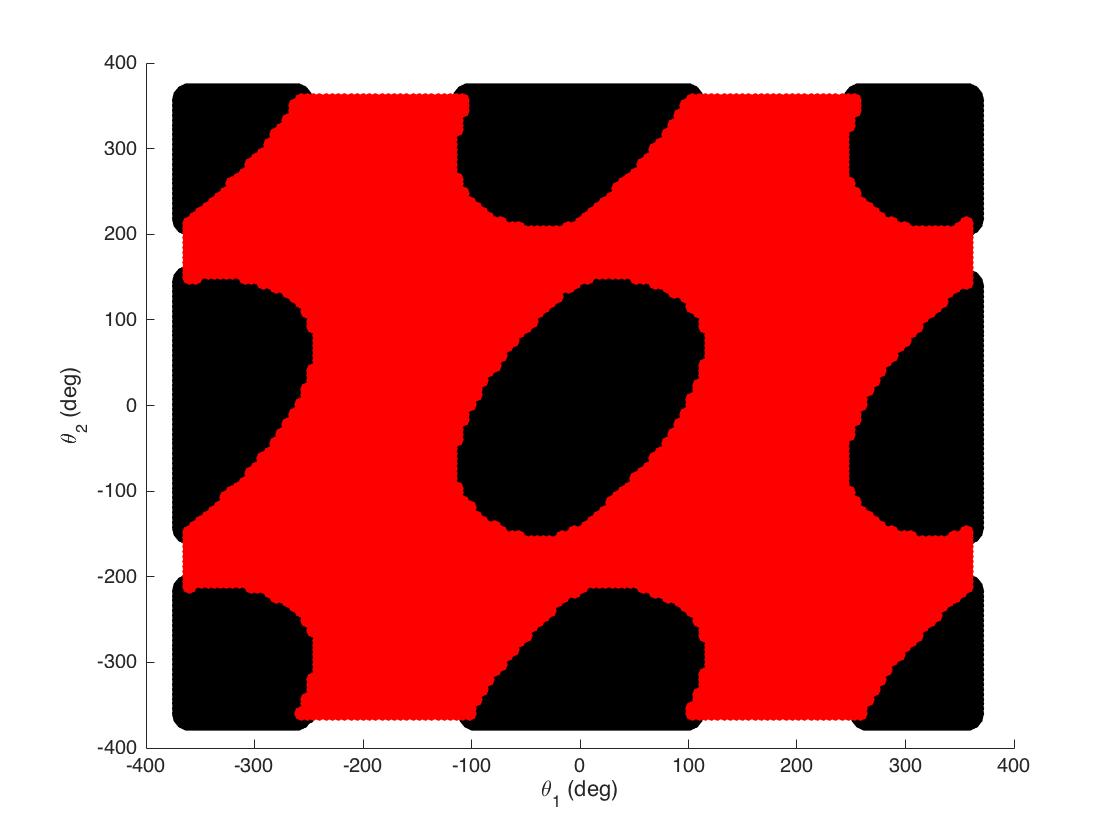}
      \caption{An inner approximation of the ROA is indicated by the back region. The projection shown is for $(\omega_1,\omega_2) = (0,0)$.}
      \label{fig:inner}
   \end{figure}


With final time $T=8$ and radius $\epsilon = 0.1$, we find the following polynomial, $v_5(0,x)$, \mbox{at fifth-order relaxation ($k=5$):} \vspace{-2em}

\begin{align*}
v_5(0,x) = & \hphantom{-}1.8707 - 4.9538x_1 + 0.0017x_2  \\
& \quad \hphantom{---------} \vdots  \\
& \quad - 0.0002x_5^2 x_6^{8} -0.0021x_5x_6^9 -0.0003x_6^{10},
\end{align*}
whose zero level set $\{ ~ x \in \mathbb{R}^6 ~|~ v_5(0,x) \geqslant 0 ~ \}$
provides an outer approximation to the ROA. We illustrate the polynomial $v_5(0,\cdot)$ in Fig.~\ref{fig:level} as a function of the original state variables $(\theta_1,\,\theta_2)$. We consider $(\omega_1,\,\omega_2) = (0,\,0)$ in order to visualize the ROA, but this is not a necessary restriction. We illustrate the outer approximation to the ROA in Fig.~\ref{fig:outer}.

Likewise, with $T=8$ and $\epsilon = 0.1$, we find at the third-order relaxation ($k=3$) the inner approximation to the ROA presented in Fig.~\ref{fig:inner} (again with $(\omega_1,\,\omega_2) = (0,\,0)$ used only for representation purposes).

We next show how one could use Hermitian SOS to obtain better numerical results. For optimal power flow problems, applying Hermitian SOS yields computational advantages while preserving convergence guarantees~\cite{josz2018}. The idea is to exploit the structure that comes from alternating current physics in order to reduce the computational burden. 
We consider the transient dynamics of a system after the fault has disappeared and we assume that there is no voltage instability. In that case, it is reasonable to assume that the magnitudes $\left|v\right|$ of the complex voltages are fixed such that only the phase angles $\theta$ are variables. This allows us to define $v_k := \exp (j \theta_k)$ (up to proper rescaling), such that $\dot{v}_k = j \dot{\theta}_k \exp (j \theta_k)$, where $j=\sqrt{-1}$. The dynamics can thus immediately be written as a differential algebraic system of equations:\vspace{-1.5em}

\begin{equation} 
\arraycolsep=0.3pt\def\arraystretch{1.5}
\left\{
\begin{array}{ccl} 
\dot{v}_k & = & \hphantom{-}j \omega_k v_k, \\
\dot{\omega}_k & = & - \lambda_k \omega_k + \frac{1}{M_k}\left(P_k - \frac{1}{2}\sum_{l \neq k} -G_{kl} |v_k|^2 - \overline{Y}_{kl} v_k \overline{v}_l - Y_{kl} v_l \overline{v}_k \right), \\
0 & = & \hphantom{-}|v_k|^2 - 1,
\end{array}
\right.
\end{equation}
where $Y_{kl}$ denotes the mutual admittance of the line connecting buses~$k$ and~$l$.

It is straightforward to adapt the theory of OMs to complex states by leveraging recent results in complex algebraic geometry \cite{angelo-2008}. Our ongoing research is implementing a complex version of the hierarchy proposed by Henrion and Korda~\cite{korda2014} in order to reduce the computational burden at a given relaxation order.



\section{Conclusion}

\label{sec:CONCLUSION}

In the context of the transient stability analysis of power systems, this paper demonstrates the potential for using the theory of occupation measures (along with convex optimization techniques) to compute inner and outer approximations to the region of attraction for a stable equilibrium point. To the best of our knowledge, this is the first time that occupation measure theory has been applied to analyze transient stability problems for electric power systems. The resulting approximations have the potential to provide analytically rigorous guarantees that can preclude the need for computationally expensive transient simulations. With computational tractability remaining an important challenge, future research will investigate how to exploit sparsity when using occupation measures.

\section*{Acknowledgment}

We wish to thank Didier Henrion at LAAS-CNRS and Milan Korda at UCSB for fruitful discussions.


\bibliography{mybib}{}
\bibliographystyle{IEEEtran}

\end{document}